\documentclass[twocolumn,10pt]{article}
\usepackage{amsmath}
\usepackage{mathtools}
\usepackage{lipsum}
\usepackage[mathscr]{euscript}
\usepackage{amsfonts}
\usepackage{amssymb}
\usepackage[none]{hyphenat}
\usepackage[margin=1in]{geometry}
\usepackage{fancyhdr}
\usepackage{graphicx}
\usepackage{layouts}
\usepackage{tikz}
\usepackage{sidecap}
\usepackage{multicol}
\usepackage{amssymb}
\usepackage[nottoc,numbib]{tocbibind}
\usepackage{mathtools}
\usepackage{float}
\usepackage{lipsum}

\usetikzlibrary{decorations.markings}
\graphicspath{{Images}}
\allowdisplaybreaks
\DeclareSymbolFont{rsfs}{U}{rsfs}{m}{n}
\DeclareSymbolFontAlphabet{\mathscrsfs}{rsfs}
\begin{document}

\begin{titlepage}
\noindent \Large{\textbf{Tensors Over Banach Spaces: 
\\ Completeness, Algebra, and Interpolation}}\\
\small{Sid Edwards}\\
\section{Abstract}
Can polynomial interpolation be extended to a Banach space setting? Are tensors whose elements are non-commutative Banach space elements legitimate objects with notable analytic and algebraic properties? Here we explore these questions and  present the ground work conducted to successfully interpolate a left-sided polynomial whose inputs are general Banach space elements, when the Banach space admits an associative algebra, submultiplicative norm and a multiplicative identity. We start with introductions and propose theorems and lemmas equipped for general tensors over a Banach space for those whose research interests lie outside of Lagrange interpolation, transitioning into algebraic properties inspired by Maksimov's work on cubic matrices, gradually narrowing down to the interpolation itself for left-sided polynomials in the last section of the article. 
\\ \\
\textbf{Key phrases}: algebra, analysis, associative,  Banach space, complete, functional, identity, induction, interpolation, inverse, linear, matrix, norm, polynomial, sequences and series, tensor, vector
\section{Acknowledgements}
The author would like to thank  Daniel Edwards for providing personal and financial support, Dr. Sigurd Angenent at UW-Madison for providing feedback on proof writing, University of Wisconsin - Madison in general for providing necessary background knowledge through several courses of undergraduate material, Dr. Rozikov Utkir at Institute of Mathematics in Tashkent, Uzbekistan for their correspondence on cubic matrices and further helpful feedback on refining this article, Gleb Pogudin at Laboratoire d'informatique de l'École polytechnique for providing constructive comments on algebraic properties, and, Chris Tisdell for providing preliminary public information related to the completeness concept on Euclidean spaces.
\tableofcontents
\vfill
\begin{center}
\line(1,0){400}\\
\small{\textbf{sidney41@usf.edu}}
\end{center}
\end{titlepage}
\newpage
\section{Space of Tensors}
\subsection{Introduction}
For general scientists who work with tensor spaces, a \textit{Banach space} is a normed vector space in which every Cauchy Sequence converges \cite{RudinF}, essentially a generalization of Euclidean vector spaces to more complex spaces of objects, such as of functions (continuously differentiable for example) or operators (matrices for example). In our case, we first generalize the concept of completeness to tensors of general Banach-valued objects. We do this for two reasons: one is because physicists may often work with tensor spaces \cite{Hess}, and two, the extension of hybercubic matrices has also been studied by other authors \cite{LaRoz} \cite{Maksimov} and as such may warrant analytical theorems compatible with their framework. 
\\ \\
An \textit{algebra} as defined in \cite{Rudin3} over a Banach space equips the Banach space with a form of multiplication such that sums and products are also elements of that space. Further, this algebra is said to be \textit{associative} \cite{Artin} if for any three elements in the space $x_1, x_2$ and $x_3,$ it is true that $(x_1x_2)x_3 = x_1(x_2x_3).$
\\ \\
To start, let $(X,|| \cdot ||_X)$ be a Banach space over the field of real numbers $\mathbb{R},$ then the tensor space $Y = X^{n_1 \times n_2 \times ... \times n_k},$ which denotes the space of tensors whose entries lie in $X,$ is also complete under a suitable norm. As our candidate, we choose the norm (for $p \geq 1$)
$$ ||y||_{Y} = \left( \sum_{\overline{i}}||y_{\overline{i}}||_{X}^{p} \right)^{1/p}.$$
It's important we take time to unpack the interpretation of these indices in order to avoid confusion. One should consider these indices a direct analog of matrix notation. Let $y$ be a rank $k$ tensor. Then the $\overline{i}$ alone in this context is not a scalar but rather a list, generally denoted as $\overline{i} = i_1, i_2, i_3 , ... ,i_k$ corresponding to entries of the rank $k$ tensor. For a specific element $y$ of $Y$, the notation $y_{\overline{i}}$ denotes any singular entry of the tensor $y$ with the commas removed from the list, $y_{\overline{i}} = y_{i_1 i_2 i_3 ... i_k }$ to form an indexed entry. Later in this document, $(y_{\overline{i}})_m$ will denote an infinite sequence in the $i_1 i_2 ... i_k$ entry of the tensor $y$. 
\\ \\
Then, $\sum_{\overline{i}}$ specifically in this document denotes multiple summations over all possible entries of the tensor as
$$\sum_{\overline{i}} = \sum_{i_1 =1}^{n_1} \sum_{i_2 = 1}^{n_2} ... \sum_{i_{k-1}=1}^{n_{k-1}} \sum_{i_k = 1}^{n_k}$$
Under these premises, we now show the above functional is indeed a norm on $Y$.
\subsection{The Norm}
Here the properties of a norm are shown to be satisfied as in \cite{RudinF}. The non-negativity of the norm is directly carried over by the norm on $X$, as we have a sum of 0. 
\\ \\
For absolute homogeneity, $\alpha \in \mathbb{R}$, we have that each entry of $y$, denoted $y_{\overline{i}}$ is an element of $X$ by construction,
\begin{align*}
||\alpha y||_{Y} &= \left(\sum_{\overline{i}}|| \alpha y_{\overline{i}}||_{X}^{p} \right)^{1/p}
\intertext{which by the norm on $X$ means}
&= \left(\sum_{\overline{i}}|\alpha|^{p}|| y_{\overline{i}}||_{X}^{p} \right)^{1/p} \\
&= \left(|\alpha|^{p} \sum_{\overline{i}}||  y_{\overline{i}}||_{X}^{p} \right)^{1/p} \\
&= |\alpha| \left( \sum_{\overline{i}}||  y_{\overline{i}}||_{X}^{p} \right)^{1/p}
\end{align*}
For the triangle inequality, we start with
\begin{align*}
||y+z||_{Y} &= \left( \sum_{\overline{i}} ||y + z||_{X}^{p} \right)^{1/p}
\intertext{which from Minkowski's inequality for non-negative numbers \cite{Mul} means}
\left( \sum_{\overline{i}} ||y + z||_{X}^{p} \right)^{1/p} &\leq \left( \sum_{\overline{i}} ||y||_{X}^{p} \right)^{1/p} + \left( \sum_{\overline{i}} ||z||_{X}^{p} \right)^{1/p} \\
\intertext{which by definition equates to}
&= ||y||_{Y} + ||z||_{Y}.
\end{align*}
\subsection{Supplemental Lemmas}
\textit{Lemma 1:} If $(y)_m$ is a Cauchy sequence in the tensor space $Y$, then each component sequence $(y_{\overline{i}})_m$ is also Cauchy.
\\ \\
\textit{Proof:} Since $(y)_m$ is Cauchy, we may pick any $\epsilon > 0$ such that for $m_1, m_2 > N$, $||(y)_{m_1} - (y)_{m_2}||_{Y} < \epsilon$. Fix any tensor entry $y_{\overline{i}}$ of the tensos $y$. Then since $y_{\overline{i}}$ is fixed, denote the sum over all tensor entries as $\overline{j} = \sum_{j_1=1}^{n_1}...\sum_{j_k=1}^{n_k}.$ We note this to suggest the norm of one fixed entry as defined originally is, separately, less than the sum of norms over all possible entries.
\begin{flalign*}
&||(y_{\overline{i}})_{m_1} - (y_{\overline{i}})_{m_2}||^{p}_{X} \leq \sum_{\overline{j}} || (y_{\overline{j}})_{m_1} - (y_{\overline{j}})_{m_2} ||_{X}^{p}
\intertext{which by the monotonicity of 1/p powers over $[0,\infty)$, $p \geq 1$ means}
&||(y_{\overline{i}})_{m_1} - (y_{\overline{i}})_{m_2}||_{X} \leq \left( \sum_{\overline{j}} ||(y_{\overline{j}})_{m_1} - (y_{\overline{j}})_{m_2}||_{X}^{p} \right)^{1/p}
\intertext{which by assumption reduces to}
&||(y_{\overline{i}})_{m_1} - (y_{\overline{i}})_{m_2}||_{X} \leq ||(y)_{m_1} - y_{m_2}||_{Y} < \epsilon.
\end{flalign*}
\\
\textit{Lemma 2:} Norms are continuous functionals.
\\ \\
\textit{Proof:} Consider a normed space $Y$ and $y_0 \in Y$, $y \in Y$ and equip $\mathbb{R}$ with the standard metric $d(a,b) = |a-b|, \ a,b, \in \mathbb{R}.$ Let $\delta > 0,$ then for $||y_0 - y||_{Y} < \delta$, we have by the reverse triangle inequality
$| ||y_0|| - ||y|| | < || y_0 - y|| < \delta.$ Now let $\delta = \epsilon,$ then we have shown norms are continuous over any $y_0 \in Y$ under the standard metric on $\mathbb{R}.$ 
\subsection{Completeness}
\textit{Theorem 1:} Let $X$ be a Banach space. Then the tensor space $Y = X^{n_1 \times n_2 \times ... \times n_k}$ is complete under the norm $$ ||y||_{Y} = \left( \sum_{\overline{i}}||y_{\overline{i}}||_{X}^{p} \right)^{1/p}.$$
\textit{Proof:} Let $(y)_m$ be a Cauchy sequence in $Y$ and let $y$ by a tensor in $Y$ whose entries are elements of $X$. By \textit{lemma 1} we know each component sequence is Cauchy in $X$, and since $X$ is complete, each component sequence converges.
\\ \\
This allows us to consider
$$\left( \sum_{\overline{i}} || (y_{\overline{i}})_m - y_{\overline{i}} ||_{X}^{p} \right)^{1/p}$$
and take the limit
$$\underset{m \rightarrow \infty}{\lim} \left( \sum_{\overline{i}} || (y_{\overline{i}})_m - y_{\overline{i}} ||_{X}^{p} \right)^{1/p}.$$ Since we are dealing with compositions and sums of continuous functions and convergent limits, we are afforded many convenient properties. By \cite{Rudin2}, the composition of continuous functions from the metric spaces $f_1: X_1 \rightarrow X_2$ and $f_2: X_2 \rightarrow X_3,$ is also continuous, and again by \cite{Rudin2}, we know $\underset{m \rightarrow \infty}{\lim}\sum_{\overline{i}}||(y_{\overline{i}})_m - y_{\overline{i}}|| =  \sum_{\overline{i}}\underset{m \rightarrow \infty}{\lim}||(y_{\overline{i}})_m - y_{\overline{i}}||.$ 
\\ \\
By this criterion, we may claim the following: if $f$ is a continuous mapping from the metric spaces $X_1$ to $X_2,$ and the sequence $(x)_m$ converges to some point $x \in X_1,$ then by \cite{Rudin2} we may claim $f(\underset{m \rightarrow \infty}{\lim}(x)_m) = f(x).$ 
\\ \\
Since each component sequences converges, and, since $p$ and $1/p$ powers are continuous functions on $[0, \infty),$ we may also leverage that this limit is equivalent to the powers of the limit points in $\mathbb{R}$ since norms are continuous mappings from a normed space to $\mathbb{R}$ by \textit{lemma 2} and compute
$$\underset{m \rightarrow \infty}{\lim} \left( \sum_{\overline{i}} || (y_{\overline{i}})_m - y_{\overline{i}} ||_{X}^{p} \right)^{1/p}$$ 
$$\left( \sum_{\overline{i}} \underset{m \rightarrow \infty}{\lim}|| (y_{\overline{i}})_m - y_{\overline{i}} ||_{X}^{p} \right)^{1/p}$$ 
$$= \left( \sum_{\overline{i}}|| \underset{m \rightarrow \infty}{\lim}( (y_{\overline{i}})_m - y_{\overline{i}} )||_{X}^{p} \right)^{1/p} = 0$$
and therefore the sequence $(y)_m$ converges to $y.$ 
\section{Algebra}
In order to establish polynomials of Banach tensors, we need to construct a notion of how to multiply them. While their multiplication in more general cases may be possible, we focus on hypercubic tensors and establish that if the space $X$ admits an associative Banach algebra with an identity, then so does the tensor space $$X^{\underbrace{n \times n \times ... \times n}_{k-times}}$$ 
by direct construction of the multiplication. 
\subsection{The Multiplication and Properties}
Much like matrices, some readers may seek an extended definition of multiplication which maps hypercubic tensors to other hypercubic tensors of the same rank. In the general case, there may be dozens upon dozens of forms of multiplication for any given hypercubic tensor, but we choose one in particular which simultaneously carries over the associativity, invertibility, submultiplicativity (to a limited extent) and a multiplicative identity.
\\ \\
Let $A$ and $B$ be two rank-$k$, $n \times n \times... \times n$ hypercubic tensors, then we define the product entry-wise. Much like matrix multiplication, we re-write components of the indices of the tensors to denote the specific entries that change and which we sum over, denoted $v$, while other indices stay the same with repsect to $v$ as they would in \cite{Artin}:
\begin{align}(AB)_{\overline{i}} = \sum_{v=1}^{n}a_{i_1 i_2 ...i_{k-1}v}b_{v i_2 ... i_k}
\end{align}
Here, $a_{i_1 i_2 i_3 ...i_{k-2}i_{k-1} v}$ contains exactly as many indices as $a_{i_1 i_2 i_3 ... i_k}$ except that we allow the last index $v$ to vary over $1 \leq v \leq n,$ independently in place of $i_k$ while $i_1 i_2 i_3...i_{k-1}$ remain fixed. Similarly, for the entry $b,$ we replace $i_1$ with $v$ and allow $v$ to vary while the remaining indices $i_2i_3...i_k$ remain fixed. The construction of this form of multiplication is inspired by one of Maksimov's studied algebras for cubic matrices in their application of stochastic processes \cite{LaRoz}\cite{Maksimov} and acts as one type of analog of matrix multiplication.
\\ \\
\textit{Theorem 2:} Let $X$ be a Banach space equipped with an associative algebra. Then the above tensor multiplication is associative over the Banach space $$X^{\underbrace{n \times n \times ... \times n}_{k-times}}.$$
\textit{Proof:} Our goal is to show for three hypercubic tensors $A,B,C$ that $(AB)C = A(BC).$ Let $a_{\overline{i}} \in \{ a_{i_1 i_2 ... v_2}: 1 \leq v \leq n \}$ and $b_{\overline{i}} \in \{b_{v_2i_2i_3...i_{k-2}i_{k-1}v_1}: 1 \leq v_1 \leq n, 1 \leq v_2 \leq n \} $ where $v_1$ and $v_2$ vary independently of the remaining indices but which take the place of the last and/or first terms of each respective tensor entry. We denote this so that we may later substitute this result into a more complex sum. Starting with the $AB$ side,
\begin{flalign}
&(AB)_{i_1...i_{k-1}v_1} = \sum_{v_2 = 1}^{n} a_{i_1i_2 ... i_{k-1}v_2}b_{v_2i_2...i_{k-1}v_1}
\intertext{where $a_{\overline{i}},b_{\overline{i}},c_{\overline{i}} \in X$ denote the respective entries of $A, B$ and $C$, then calculate}
&((AB)C)_{\overline{i}} = \sum_{v_1=1}^{n}(AB)_{i_1...i_{k-1}v_1}c_{v_1i_2...i_{k-1}i_k}
\intertext{and substitute the result from $(2)$ into $(3)$  above}
&((AB)C)_{\overline{i}} = \\
\nonumber &\sum_{v_1=1}^{n}\left( \sum_{v_2=1}^{n}a_{i_1...i_{k-1}v_2}b_{v_2i_2...-{k-1}v_1} \right)c_{v_1i_2...i_{k-1}i_k}
\intertext{which by the associativity inherited by $X$ gives us}
&= \sum_{v_1=1}^{n}\sum_{v_2=1}^{n}a_{i_1i_2...i_{k-1}v_2}b_{v_2i_2...i_{k-1}v_1}c_{v_1i_2...i_{k-1}i_k}
\end{flalign}
Now we approach the other side:
\begin{flalign}
&(BC)_{v_2i_2...i_{l-1}i_k} = \sum_{v_1=1}^{n}b_{v_2i_2...i_{k-1}v_1}c_{v_1i_2...i_k} 
\intertext{then define}
&A(BC)_{\overline{i}} = \sum_{v_1=1}^{n}a_{i_1...v_2}(BC)_{v_2i_2...i_k}
\intertext{and then substitute $(5)$ into $(6)$}
&A(BC)_{\overline{i}} = \\
\nonumber & \sum_{v_1=1}^{n}a_{i_1i_2...i_{k-1}v_2}\left( \sum_{v_1=1}^{n}b_{v_2i_2...i_{k-1}v_1}c_{v_1i_2...i_k} \right)
\intertext{where by the associativity of $X$, we find $(7)$ simplifies to}
&= \sum_{v_1=1}^{n}\sum_{v_1=1}^{n}a_{i_1...i_{k-1}v_2}b_{v_2i_2...i_{k-1}v_1}c_{v_1i_2...i_k}
\end{flalign}
then we have shown $(AB)C = A(BC).$
\\ \\
\textit{Theorem 3:} Suppose the Banach space X admits an associative algebra with an identity I, then the space
$$Y = X^{\underbrace{n \times n \times ... \times n}_{k-times}}.$$
also admits a Banach algebra with an identity. 
\\ \\
\textit{Proof:} We are tasked with finding a compatible identity element for any rank of hypercubic tensor. 
\\ \\
Define a tensor $T \in Y$ with the properties
$$T_{\overline{i}} = \begin{cases} I & \text{if} \ i_1 = i_k \\ 0 & \text{if} \ i_1 \neq i_k \end{cases}$$ where $I$ is the multiplicative identity of the Banach space $X$ and $0$ is the additive identity on $X.$
\\ \\
Now take any other tensor $y \in Y$ and multiply it with $T$ to obtain
\begin{align*}
(yT)_{\overline{i}} &= \sum_{v=1}^{n} y_{i_1...i_{k-1}v}T_{vi_2...i_k}
\intertext{Recall that $v$ may vary while other indices remain fixed, then if we truncate the sum (where $v$ takes the place of $i_k$ or $i_1$)}
&= \sum_{v = 1}^{i_{k}-1}y_{i_1i_2...i_{k-1}v}T_{vi_2i_3...i_k} \\
&+ \sum_{v=i_k}^{i_k}y_{i_1i_2...i_{k-1}v}T_{vi_2...i_k} \\
& + \sum_{v=i_k + 1}^{n}y_{i_1i_2...i_{k-1}v}T_{vi_1i_2...i_k}
\intertext{we find by construction and the uniqueness of each tensor entry that the above equates to}
& =\sum_{v=1}^{i_k - 1}y_{i_1i_2...i_{k-1}v} \cdot 0 + y_{\overline{i}} \cdot I \\
&+ \sum_{v = i_{k}+1}^{n}y_{i_1...i_{k-1}v} \cdot 0 \\
& = y_{\overline{i}}
\end{align*}
A nearly identical process shows the left-side identity with $Ty = y$ by truncating the sum at $v = i_1$ instead. This work shows the construction and existence of an identity element for any tesnsor space over a Banach space $Y.$ The uniqueness of this identity can be established as follows: suppose $T_1$ and $T_2$ are two identities of the hypercubic tensor space $Y$. Then  $T_1 = T_1 T_2$ and $T_1 T_2 = T_2 \implies T_1 = T_2.$ 
\subsection{Submultiplicativity}
Now that we have multiplication defined, another goal is to prove that if the Banach space $X$ has a submultiplicative norm, then there exists a submultiplicative $2$-norm on the hypercubic tensor space $Y = X^{n \times n \times ... \times n}.$
\\ \\
A norm $||\cdot||$ is submultiplicative on the Banach space $X$ if, for two elements $x_1$ and $x_2$, it is true that $||x_1 x_2||_X \leq ||x_1||_X||x_2||_X$ \cite{Lan}. This type of norm is particularly useful in the study of matrices and will be important later on when we talk of inversion. It is trivial to show there exists at least one submultiplicative norm on $Y$ when we take $p = 1$, but we go beyond that to accomodate more numerical methods and approximations which rely on squares \cite{Haibn}. We will use two inductive proofs to accomodate theorem 4 to conclude the final result.
\\ \\
\textit{Theorem 4:} Let $X$ be a Banach space with an associative algebra equipped with a submultiplicative norm $|| \cdot ||_X$ and identity $I.$ Then the $2$-norm of the hypercube tensor space $Y$ over $X$, denoted $||\cdot||_{Y,2}$ with $p=2$ is also submultiplicative.
\\ \\
\textit{Proof:} Let $A_j, B_j \in X.$ Our first inductive argument will be to show that for a single sum with some singular index $j$ with a defined multiplication between elements yields
$$\sum_{j=1}^{n}||A_jB_j||^2 \leq \left( \sum_{j=1}^{n}||A_j||^2 \right) \cdot \left( \sum_{j=1}^{n}||B_j||^2 \right).$$
As our base case for $n=1$ we have by the montonicity of squares over $[0,\infty]$ that $||A_jB_j||^2 \leq ||A_j||^2||B_j||^2.$ Now assume the case holds for any integer $m$, that is
$$\sum_{j_=1}^{m}||A_jB_j||^2 \leq \left( \sum_{j=1}^{m}||A_j||^2 \right) \cdot \left( \sum_{j=1}^{m}||B_j||^2 \right),$$
then we take
\begin{flalign*}
& \sum_{j=1}^{m+1}||A_jB_j||^2 = \sum_{j=1}^{m}||A_jB_j||^2 + ||A_{m+1}B_{m+1}||^2 
\intertext{and find by assumption that}
& \sum_{j=1}^{m}||A_jB_j||^2 + ||A_{m+1}B_{m+1}||^2 \leq  \\
& \left( \sum_{j=1}^{m}||A_j||^2 \right) \cdot \left( \sum_{j=1}^{m}||B_j||^2 \right) + ||A_{m+1}B_{m+1}||^2 
\intertext{Next we complete the square of the $m+1$}
& \left( \sum_{j=1}^{m}||A_j||^2 \right) \cdot \left( \sum_{j=1}^{m}||B_j||^2 \right) + ||A_{m+1}B_{m+1}||^2 \leq  \\
& \left( \sum_{j=1}^{m}||A_j||^2 \right) \cdot \left( \sum_{j=1}^{m}||B_j||^2 \right) + \\
& ||A_{m+1}||^2 \left(\sum_{j=1}^{m}||B_j||^2 \right) +||B_{m+1}||^2 \left( \sum_{j=1}^{n}||A_j||^2 \right) \\
& + ||A_{m+1}||^2||B_{m+1}||^2
\intertext{which equates to}
& = \left( \sum_{j=1}^{m+1}||A_j||^2 \right) \cdot \left( \sum_{j=1}^{m+1}||B_j||^2 \right).
\end{flalign*}
Now let $m = n$, then we are done with this inductive step.
\\ \\
For our next inductive argument, we extend the preceding results to a more general tensor structure. We next prove (where $\overline{i}$ here denotes the same list of tensor entries as defined in the first section of this article)
$$\sum_{i_1}^{n_1}...\sum_{i_k}^{n_k}||A_{\overline{i}}B_{\overline{i}}||^2 \leq \left( \sum_{\overline{i}}||A_{\overline{i}}||^2 \right) \left( \sum_{\overline{i}}||B_{\overline{i}}||^2 \right)$$
For our base case, we now know 
$$\sum_{i_1 = 1}^{n_1}||A_{i_1}B_{i_1}||^2 \leq \left( \sum_{i_1=1}^{n_1}||A_{i_1}||^2 \right) \left( \sum_{i_1 = 1}^{n_1} ||B_{i_1}||^2 \right)$$
Now assume the $m$-th case holds and consider
\begin{flalign*}
& \sum_{\overline{i}} \sum_{i_{m+1}}^{n_{m+1}}||A_{(i)i_{m+1}}B_{(i)i_{m+1}}||^2
\intertext{where $\sum_{\overline{i}} =  \sum_{i_1 =1}^{n_1} ... \sum_{i_{m}=1}^{n_m}.$}
\intertext{Then by our assumption of the submultiplicativity of the norm on $X,$ the above is less that or equal to}
&  \leq \sum_{\overline{i}}\left( \left( \sum_{i_{m+1}}^{n_{m+1}} ||A_{(i)i_{m+1}}||^2 \right)\left( \sum_{i_{m+1}}^{n_{m+1}} ||B_{(i)i_{m+1}}||^2 \right) \right)
\intertext{which is further}
& \leq \left( \sum_{\overline{i}} \sum_{i_{m+1}}^{n_{m+1}}||A_{(i)i_{m+1}}||^2 \right) \left( \sum_{\overline{i}} \sum_{i_{m+1}}^{n_{m+1}}||A_{(i)i_{m+1}}||^2 \right).
\end{flalign*}
Let $n_m$ = $n_k,$ the inductive argument is complete.
\\ \\
With these two inductive components proven, we can finish \textit{theorem 4} as follows:
\begin{flalign*}
& ||A_{\overline{i}}B_{\overline{i}}|| \leq ||A_{\overline{i}}||||B_{\overline{i}}|| \\
& ||A_{\overline{i}}B_{\overline{i}}||^2 \leq ||A_{\overline{i}}||^2 ||B_{\overline{i}}||^2 \\
& \sum_{\overline{i}} ||A_{\overline{i}}B_{\overline{i}}||^2 \leq \sum_{\overline{i}} ||A_{\overline{i}}||^2 ||B_{\overline{i}}||^2 \\
& \leq \left(\sum_{\overline{i}}||A_{\overline{i}}||^2 \right) \left( \sum_{\overline{i}}||B_i||^2 \right)
\intertext{then we take the square root over the inequality}
& \left( \sum_{\overline{i}} ||A_{\overline{i}}B_{\overline{i}}||^2 \right)^{1/2} \leq \\
& \left(\sum_{\overline{i}}||A_{\overline{i}}||^2 \right)^{1/2} \left( \sum_{\overline{i}}||B_i||^2 \right)^{1/2}
\end{flalign*}
showing $||AB||_{Y,2} \leq ||A||_{Y,2} ||B||_{Y,2}.$
\subsection{Inverse}
Here we accomodate the existence of inverses in such a general setting as tensors of non-commutative Banach elements by discussing a local, analytic method. If $y \in Y$ where $Y$ is a hypercube tensor space over a submultiplicative Banach space $X$ equipped with an associative algebra and identity, then one may define the point-wise convergent inverse
$$y^{-1} = \sum_{j=0}^{\infty} (-1)^{j}(y - I)^j $$ 
where $I$ in this context is the multiplicative identity of the tensor space $Y$ whose existence was proven earlier in this document using the tensor $T.$ Since we may pick a submultiplicative norm on $Y$, it is trivial to show this series converges when $||y-I||_Y < 1$ via submultiplicativity and the triangle inequality. However, showing that this series possess the desired functional property that $y \cdot y^{-1} = I$ is another matter that remains to be proven. To prove this series functions as an inverse, we take inspiration from Wu et al.'s work with the Neumann series \cite{Wu}.
\\ \\
\textit{Theorem 5:} Let $Y$ be a tensor space over a Banach space $X$ with an associative algebra, submultiplicative norm and identity. Define a series $S = \sum_{j=0}^{\infty} (-1)^{j}(y - I)^j$ with $||y - I||_Y < 1.$ Then $yS = I.$
\\ \\
\textit{Proof:} We first start with showing that the action of this series on an element $y - I$ converges by considering the truncated series $S_m.$ Consider a submultiplicative norm on $Y$, then
\begin{flalign}
\nonumber & ||(y-I)S_m - (y-I)S||_Y \leq ||y-I||_Y \cdot ||S_m - S||_Y
\intertext{where by the assumption of $||y-I||_Y<1,$}
& ||(y-I)S_m - (y-I)S||_Y \leq ||S_m - S||_Y.
\intertext{taking the convergence of the original series yields}
& \underset{m \rightarrow \infty}{\lim}||S_m - S||_Y = 0
\intertext{and thus}
& \underset{m \rightarrow \infty}{\lim}||(y-I)S_m - (y-I)S||_Y \leq 0.
\end{flalign}
Now we know $(y-I)S$ is a conergent series. It would be helpful however if we also knew that $(y-I)S$ converged to another algebraic structure, specifically that $(A-I)S = I-S$. For this, we start again with a truncated series and take $(y-I)^{0} \equiv I,$
\begin{flalign}
&(y - I) \cdot \sum_{j=1}^{m}(-1)^{j}(y-I)^{j} = \sum_{j=0}^{m}(-1)^{j}(y-I)^{m+1}
\intertext{and relabel $j+1 = l$ which implies $j= l-1,$ yielding}
&\sum_{l=1}^{m+1}(-1)^{l-1}(y-I)^{l} + I - I = \\
\nonumber & -\sum_{l=0}^{m+1}(-1)^{l}(y-I)^l + I
\intertext{which equates to}
&I - S_{m+1} \implies (y-I)S_m = I - S_{m+1}.
\intertext{We now combine results from $(15)$ and $(12)$:}
& \underset{m \rightarrow \infty}{\lim}(y-I)S_m = \underset{m \rightarrow \infty}{\lim} (I - S_{m+1})
\intertext{where since the series $S_{m}$ converges, we then know}
& \underset{m \rightarrow \infty}{\lim}( I - S_{m+1}) = I - S
\intertext{and therefore}
& (y-I)S = I - S.
\intertext{from here,}
\nonumber & (y-I)S + S = I \\
\nonumber & (y-I + I)S = I \\
& yS = I.
\end{flalign}
\section{Lagrange Interpolation}
With the work of general Banach tensors accounted for, we can narrow down to the cases of Banach matrices and their action on Banach column vectors to generalize Lagrange's method of polynomial interpolation. 
\\ \\
\textit{Theorem 6:} Let $X$ be a Banach space with an assocative algebra, submultiplicative norm and identity. Then multiplication of a Banach column vector $x \in X^{n}$ by an invertible rank $k$ hypercube tensor $M \in X^{n \times n \times ... \times n} = Y$ is a bijective linear map from $X^{n} \rightarrow X^{n}.$ 
\\ \\
\textit{Proof:} 
We prove bijectivity by proving injectivity and surjectivity assuming a definition of multiplication that will be explained later in this proof. For injectivity, take any $x \in X^n$ and $M \in Y$ and assume $Mx = Mz,$ then we may take $M^{-1}Tx = M^{-1}Mz \implies x = z.$
\\ \\
For surjectivity, consider any $x \in X^{n},$ then since $M$ is invertible, we have $MM^{-1}x = x.$ Now let $z = M^{-1}x,$ then we have $Mz = x.$ Since $Mz \in X^n$ and $M$ is a map from $X^n \rightarrow X^n,$ it follows there exists $z \in X^{n}$ such that $Mz = x.$ 
\\ \\
On linearity, let $x, z \in X^{n},$ $a, b \in \mathbb{R}$  and let $M \in Y,$ then our goal is to show $M(ax + bz) = aM(x) + bM(z).$
\\ \\
To define multiplication of a hypercube tensor on a column vector, we need only to omit the extra indices in preceding constructions and define the multiplication component-wise. Consider any tensor entry, then the components of the sum are
$$M(ax+bz)_{i} = \sum_{v=1}^{n}M_{i_1i_2...i_{k-1}v}(ax+bz)_v.$$
where here, the single scalar entry $i$ takes the place of $\overline{i}$ to denote that there is only a single index to accomodate the rank $1$ column vector. By definition, $X$ is a vector space, and by construction, each entry of the tensor $M$ and column vectors $x, \ z$ are entries of $X$ which implies they each inherit linearity. Thus,
\begin{align}
\sum_{v=1}^{n}M_{i_1...i_{k-1}v}(ax+bz)_v =\\
 \sum_{v=1}^{n}M_{i_1...i_{k-1}v}(ax)_v + \sum_{v=1}^{n}M_{i_1...i_{k-1}v}(bz)_v \\
= a(Tx)_{i} + bT(z)_{i}
\end{align}
Since this is true for any component of the vectors $x,z$ for $1 \leq i_k \leq n$, it follows that $M(ax + bz) = aM(x) + bM(z).$
\\ \\
The linearity of the inverse also follows from this result: let 
\begin{flalign*}
& w = \alpha M^{-1}x + M^{-1}z, \ \text{then} \\
& Tw = M(\alpha M^{-1}x + \beta M^{-1}z) \\
\intertext{where by the linearity of M we have}
& Tw = \alpha MM^{-1}x + \beta MM^{-1}z \\
& = \alpha x + \beta z \\
& M^{-1}Tw = w = M^{-1}(\alpha x + \beta z).
\end{flalign*}
Now we can move forward with Lagrange interpolation with these spaces. Let $x \in X$ where $X$ is a Banach space with an associative algebra, submultiplicative norm and identity. Then construct the polynomial $P(x) = \sum_{j=0}^{n} \alpha_{j}x^{j}$ with $x^{0} \equiv I,$ $\alpha_0,...,\alpha_n \in X.$ Pick $n$ data-points, $n \in \mathbb{N},$ $(x_0,P_0),(x_1,P_1)...(x_{n},P_{n})$ to construct the system of equations
\begin{flalign*}
& P_0 =  \sum_{j=0}^{n} \alpha_{j}x_1^{j} \\
& P_1 =  \sum_{j=0}^{n} \alpha_{j}x_2^{j} \\
& \ \ \ \ \ \ \ \ \ \vdots \\
& P_n =  \sum_{j=0}^{n} \alpha_{j}x_n^{j}
\end{flalign*}
which in vector-matrix form is
$$ \begin{bmatrix} P_0 \\ P_1 \\ \vdots \\ P_n \end{bmatrix} = \begin{bmatrix} 1 & x_0 & ... & x_0^n \\ 1 & x_0 & ... & x_1^n \\ \vdots & \vdots & \ddots &\vdots \\ 1 & x_n & ... & x_n^n \end{bmatrix}  \begin{bmatrix} \alpha_0 \\ \alpha_1 \\ \vdots \\ \alpha_n \end{bmatrix}$$
and represent this as $P = x\alpha.$ For any chosen vector $P$, if the matrix $x$ is invertible, then this system has a unique solution given by $x^{-1}P=\alpha.$


\begin{thebibliography}{1}
\bibitem{Haibn} H. Chen, Y. Wang, G. Zhou, \textit{High-order sum-of-squares structured tensors: theory and applications}, Frontiers of Mathematics in China, \textbf{15} (2020), no. 4, pp. 255-284. DOI: \text{https://doi.org/10.1007/s11464-020-0833-1}. 
\bibitem{Mul} H. P. Mulholland, \textit{On Generalizations of Minkowski's Inequality in the Form of a Triangle Inequality}, Proceedings of the London Mathematical Society, s2-\textbf{51} (1949), no. 1, DOI: \text{https://doi.org/10.1112/plms/s2-51.4.294}
\bibitem{Lan} P. Lancaster, M. Tismenetsky, \textit{The Theory of Matrices: With Applications}, Academic Press, United States, (1985), no. 2, pp. 358. ISBN: 0-12-435560-9.
\bibitem{Wu}
M. Wu et al., \textit{Approximate Matrix Inversion for High-Throughput Data Detection in the Large-Scale MIMO Uplink}, IEEE International Symposium on Circuits and Systems (ISCAS), 19-23 May 2013, DOI: 10.1109/ISCAS21123.2013.
\bibitem{Artin}
M. Artin, \textit{Algebra}, Pearson College Div, India, 2015, no. 1, pp. 4-5, ISBN-10: 0130047635.
\bibitem{LaRoz}
M. Ladra, U. A. Rozikov, Algebras of cubic matrices, {\it Linear and Multilinear Alg.} 2017. V.65, No.7, p. 1316--1328.
\bibitem{Hess}
S. Hess, \textit{Tensors for Physicists}, Springer, 2015, pp. 43-54, ISBN-10: 3319127861.
\bibitem{Maksimov}
V. M. Maksimov, \textit{CUBIC STOCHASTIC MATRICES AND THEIR PROBABILITY INTERPRETATION*} (Translated by M. V. Khatuntseva), Theory of Probability and its Applications, \textbf{41} (1997), no. 1, DOI: \text{https://doi.org/10.1137/TPRBAU00004100000}\\ \text{1000055000001}.
\bibitem{RudinF}
W. Rudin, \textit{Functional Analysis}, McGraw-Hill Inc, United States, 1991, no. 2, pp 4. ISBN: 0-07-054236-8.
\bibitem{Rudin}
W. Rudin, \textit{Principles of Mathematical Analysis}, McGraw-Hill Inc, United States, 1976, no. 3, pp 220. ISBN: 0-07-054235-X.
\bibitem{Rudin2}
W. Rudin, \textit{Principles of Mathematical Analysis}, McGraw-Hill Inc, United States, 1976, no. 3, pp 84-87. ISBN: 0-07-054235-X.
\bibitem{Rudin3}
W. Rudin, \textit{Principles of Mathematical Analysis}, McGraw-Hill Inc, United States, 1976, no. 3, pp 161. ISBN: 0-07-054235-X.
\end{thebibliography}
\end{document}